\documentclass{article}

\def\sym{{\mathrm{Sym} }}
\newtheorem{lem}{Lemma}
\newtheorem{cor}{Corollary}
\newtheorem{prop}{Proposition}
\def\bz{\ensuremath{\mathbf{Z}}}
\def\bc{\ensuremath{\mathbf{C}}}
\def\bp{\ensuremath{\mathbf{P}}}
\def\a{\mathbf{a}}

\title{Nondegeneracy for Quotient Varieties under Finite Group Actions}
\author{ S. S. Kannan  and P. Vanchinathan}

\begin{document}
\maketitle
\begin{center}
\large {Chennai Mathemtical Institute\\
 H1, SIPCOT IT Park \\
Padur PO, Siruseri 603103\\
India\\
email:\quad \tt kannan@cmi.ac.in, vanchi@cmi.ac.in\\ }
\end{center}
\begin{abstract}
\noindent We prove that for an abelian group $G$ of order $n$ 
the morphism 
$ \varphi\colon \bp(V^*)\longrightarrow 
\bp ( (\sym^n V^*)^G )$
defined by $\varphi( [f] ) = [\prod_{\sigma\in G} \sigma \cdot f ] $
is nondegenerate for every finite-dimensional representation
$V$ of $G$
 if and only if either $n $ is a prime number or $n=4$.

\end{abstract}


\vspace{6pt}
\textbf{AMS Classification No: 14L30}

\section*{Introduction}
Let $V$ be a finite-dimensional representation of a finite group
$G$ of order $n$ over a field $k$.  Then by a classical theorem of 
Hermann Weyl \cite{weyl} on polarizations  one obtains a nice generating set for
the invariant ring $(\sym\, V^*)^G$ consisting of 
$ e_r(f), r=1,2,\ldots,n,\  f\in V^*$ where $e_r(f)$
denotes the $r$-th elementary symmetric polynomial in
$\sigma.f,\ \sigma\in G$.

Also, for any point $x\in\bp(V)$ there is a $G$-invariant section
of the line bundle $\mathcal{O} (1)^{\otimes n}$ of the form
$s = e_n(f) = \prod_{\sigma\in G}\sigma . f ,\ f\in V^*$
such that $s(x)\neq0$. So, the line bundle $\mathcal{O}(1)^{\otimes n}$
descends to the quotient $G\backslash \bp(V)$.  (See \cite{mum}, \cite{news}).

This leads one to ask the natural question whether 
the set $\{\prod_{\sigma\in G} \sigma.f : f\in V^*\}$ generates
the $k$-algebra $\bigoplus_{d\in \bz_{\geq0} } 
  \big(\sym^{dn} V^*\big)^G$ of $G$-invariants.

In particular in view of Theorem 3.1   of \cite{kannan}  
the following question arises: Is the morphism
$$\bp(V^*)\longrightarrow \bp \big( (\sym^{n} V^*)^G\big)$$ given
by $  [f] \mapsto [\prod_{\sigma\in G} \sigma \cdot f ] $ 
nondegenerate?

\vspace{6pt}
In an attempt to answer this question we prove the following
result:

\vspace{3pt}
\noindent{\bf Theorem}:\quad{\it
Let $G$ be a finite abelian group of order $n$. Then
the map $$ \varphi\colon \bp(V^*)\longrightarrow 
\bp ( (\sym^n V^*)^G )$$
defined by $\varphi( [f] ) = [\prod_{\sigma\in G} \sigma \cdot f ] $ 
is nondegenerate
for every finite-dimensional complex representation $V$ of $G$
 if and only if $n $ is either a prime number or $n=4$.}

The paper is organized as follows:

In Section 1 some preliminaries are established.

In Section 2 we give a criterion linking nondegeneracy to
nonsingularity of an associated matrix. 

In Section 3 the main theorem is proved.
\section{Preliminaries}

Recall that a map into a vector space $V$ or $\bp(V)$ is said to be
degenerate, if the image is contained in a hyperplane.

\begin{lem}
Let $V$ be a finite-dimensional vector space
over an infinite field $k$, and $f_1,f_2,\ldots,f_N$ be a 
finite number of linearly independent elements of 
$\sym(V^*) = \bigoplus
_{d=0}^\infty \sym^d V^*$. Then the morphism $\psi\colon V\to k^N$
given by $\psi(v) = (f_1(v), f_2(v),\ldots, f_N(v) )$ is 
nondegenerate.
\end{lem}
Proof:\quad If $\psi$ were degenerate, then there would exist a nonzero
linear form $F =\sum_i \alpha_iY_i$  such that  $F(\psi(V))=0$.
Then $\sum_i\alpha_i f_i(v) = 0$ for all $v\in V$.
This implies, as the base field is infinite,
that $\sum_i\alpha_if_i=0$
contradicting the linear independence of $f_i$'s. 

\begin{lem} Let $v_1$ be a nonzero vector
in a vector space $V$ over $\bz/2\bz$ of finite-dimension. 
Then the two sets,
$\{ A\subset V \mid |A| \hbox{ is even and } 
\sum_{w\in A} w = 0\}$ and 
$\{ A\subset V \mid |A| \hbox{ is odd and } \sum_{w\in A} w = v_1\}$ 
have same cardinality.
\end{lem}

Proof:\quad The correspondence $f(A) = A\setminus\{v_1\} \hbox { or }
A\cup \{v_1\}$ according as whether $v_1$ is in $A$ or not, sets up
the bijection between these two collections of subsets of $V$.

\section{ Nondegeneracy and Nonsingularity}

First we fix some notations.

Let $G$ be a finite abelian group of order $n$, and $V$ a finite-dimensional  permutation representation for $G$. If $\{X_i\}_{i=1}^m$ is a basis of
of $V^*$ permuted by the action of $G$ then the monomials on 
$X_i$'s of degree $d$  provide a basis for  
$\sym^d V^*$ again permuted by $G$ action.

Let $O_1,O_2,\ldots, O_N$ be all the distinct $G$-orbits
for monomials of degree  $n$ (${}=|G|)$ in these $X_i$'s.
For an $r\in \{1,2,\ldots, N\}$ we 
denote by $O_r(X)$ the element in $(\sym^n V)^G$ that is
the sum of the monomials  in the orbit $O_r$.
Working with the $G$-module $V^*\oplus V^*$, we use 
the notation $Y_i$  for the  basis of the other copy of 
 $V^*$ and define $O_r(Y)$ similarly.

Now we express the $G$-invariant  
$\prod_{\sigma\in G} \sigma (\sum_{i=1}^n X_iY_i)$ as 
linear combination of these orbit sums 
$\sum_{r,s=1}^N a_{rs} O_r(X) O_s(Y)$.

\begin{lem}
A necessary and sufficient 
condition for the morphism $\varphi\colon V^*\to
(\sym ^n V^*)^G$ given by
$\varphi(f) = \prod_{\sigma\in G}  \sigma.f$ to be
 nondegenerate is that the $N\times N$-matrix $((a_{rs}))$ 
be nonsingular.
\end{lem}

Proof:\quad Assume that the matrix $((a_{rs}))$ is non-singular.
By previous lemma the map $V\to k^N$ sending $v$ to 
$(O_1(Y)(v), \ldots, O_N(Y)(v) )$ is nondegenerate.

Hence, there exists $v_1,v_2,\ldots, v_N\in V$ such
that the matrix $((O_r(Y)(v_s))$ is nonsingular. Hence the product
of the two matrices $ ((a_{rs}))\cdot ((O_r(Y)(v_s))$ is also
nonsingular.  Hence the image of the map
$V^*\to (\sym^n V^*)^G$ defined by 
$f\mapsto \prod_{\sigma\in G} \sigma.f$
contains a basis of $(\sym^n V^*)^G$.

Converse is proved similarly.

\noindent{\bf Example} For any finite group $G$ of order ${}\leq5$, and
$V$ its regular representation, the matrix $((a_{rs}))$ has determinant
$\pm1$.   

For the case of group of order 5 the regular representation
 has 26 orbits  on the set of monomials of degree 5
 on 5 variables  and  the $26\times 26$ matrix is 
as given below. 

{\footnotesize
$$\arraycolsep=4pt
\left[
\begin{array}{cccccccccccccccccccccccccc}
 0& 0& 0& 0& 0& 0& 0& 0& 0& 0& 0& 0& 0& 0& 0& 0& 0& 0& 0& 0& 0& 0& 0& 0& 0& 1 \\
 0& 0& 0& 0& 0& 0& 0& 0& 0& 0& 0& 0& 0& 0& 0& 0& 0& 0& 0& 0& 0& 1& 0& 0& 0& 0 \\
 0& 0& 0& 0& 0& 0& 0& 0& 0& 0& 0& 0& 0& 0& 0& 0& 0& 0& 0& 0& 0& 0& 1& 0& 0& 0 \\
 0& 0& 0& 0& 0& 0& 0& 0& 0& 0& 0& 0& 0& 0& 0& 0& 0& 0& 0& 0& 0& 0& 0& 1& 0& 0 \\
 0& 0& 0& 0& 0& 0& 0& 0& 0& 0& 0& 0& 0& 0& 0& 0& 0& 0& 0& 0& 0& 0& 0& 0& 1& 0 \\
 0& 0& 0& 0& 0& 0& 0& 0& 0& 0& 0& 0& 0& 0& 0& 0& 0& 0& 0& 1& 0& 0& 1& 0& 0& 0 \\
 0& 0& 0& 0& 0& 0& 0& 0& 0& 0& 0& 0& 0& 0& 0& 1& 0& 0& 0& 0& 0& 0& 0& 0& 1& 0 \\
 0& 0& 0& 0& 0& 0& 0& 0& 0& 0& 0& 0& 0& 0& 0& 0& 0& 1& 0& 0& 0& 1& 0& 0& 0& 0 \\
 0& 0& 0& 0& 0& 0& 0& 0& 0& 0& 0& 0& 0& 0& 0& 0& 0& 0& 0& 0& 1& 0& 0& 1& 0& 0 \\
 0& 0& 0& 0& 0& 0& 0& 0& 0& 1& 0& 0& 0& 0& 0& 0& 0& 0& 0& 0& 1& 0& 0& 2& 0& 0 \\
 0& 0& 0& 0& 0& 0& 0& 0& 0& 0& 1& 0& 0& 0& 0& 1& 0& 0& 0& 0& 0& 0& 0& 0& 2& 0 \\
 0& 0& 0& 0& 0& 0& 0& 0& 0& 0& 0& 0& 1& 0& 0& 0& 1& 0& 1& 0& 0& 0& 0& 0& 0& 5 \\
 0& 0& 0& 0& 0& 0& 0& 0& 0& 0& 0& 1& 0& 0& 0& 0& 1& 0& 1& 0& 0& 0& 0& 0& 0& 5 \\
 0& 0& 0& 0& 0& 0& 0& 0& 0& 0& 0& 0& 0& 1& 0& 0& 0& 1& 0& 0& 0& 2& 0& 0& 0& 0 \\
 0& 0& 0& 0& 0& 0& 0& 0& 0& 0& 0& 0& 0& 0& 1& 0& 0& 0& 0& 1& 0& 0& 2& 0& 0& 0 \\
 0& 0& 0& 0& 0& 0& 1& 0& 0& 0& 1& 0& 0& 0& 0& 0& 0& 0& 0& 0& 0& 0& 0& 0& 3& 0 \\
 0& 0& 0& 0& 0& 0& 0& 0& 0& 0& 0& 1& 1& 0& 0& 0& 2& 0& 1& 0& 0& 0& 0& 0& 0& 5 \\
 0& 0& 0& 0& 0& 0& 0& 1& 0& 0& 0& 0& 0& 1& 0& 0& 0& 1& 0& 0& 0& 3& 0& 0& 0& 0 \\
 0& 0& 0& 0& 0& 0& 0& 0& 0& 0& 0& 1& 1& 0& 0& 0& 1& 0& 2& 0& 0& 0& 0& 0& 0& 5 \\
 0& 0& 0& 0& 0& 1& 0& 0& 0& 0& 0& 0& 0& 0& 1& 0& 0& 0& 0& 1& 0& 0& 3& 0& 0& 0 \\
 0& 0& 0& 0& 0& 0& 0& 0& 1& 1& 0& 0& 0& 0& 0& 0& 0& 0& 0& 0& 1& 0& 0& 3& 0& 0 \\
 0& 1& 0& 0& 0& 0& 0& 1& 0& 0& 0& 0& 0& 2& 0& 0& 0& 3& 0& 0& 0& 5& 0& 0& 0& 0 \\
 0& 0& 1& 0& 0& 1& 0& 0& 0& 0& 0& 0& 0& 0& 2& 0& 0& 0& 0& 3& 0& 0& 5& 0& 0& 0 \\
 0& 0& 0& 1& 0& 0& 0& 0& 1& 2& 0& 0& 0& 0& 0& 0& 0& 0& 0& 0& 3& 0& 0& 5& 0& 0 \\
 0& 0& 0& 0& 1& 0& 1& 0& 0& 0& 2& 0& 0& 0& 0& 3& 0& 0& 0& 0& 0& 0& 0& 0& 5& 0 \\
 1& 0& 0& 0& 0& 0& 0& 0& 0& 0& 0& 5& 5& 0& 0& 0& 5& 0& 5& 0& 0& 0& 0& 0& 0& 15 \\
\end{array}\right]$$
}

The determinant being a unit in $\bz$
the nondegeneracy of the map $\varphi$ holds
in all characteristics.

\noindent{\bf Remark}\quad In the case of $V$ having a basis
of common eigenvectors for the action of $G$, a variation of
this lemma can be proved easily. In such cases the matrix
$((a_{rs}))$ would be diagonal. In the next section this variation
will be applied.

\section{Criterion for Nondegeneracy of Finite Abelian Group Quotients}
In this section we prove the following theorem.

\noindent{\bf Theorem}:\quad{\it
Let $G$ be a finite abelian group of order $n$. Then
the map $$ \varphi\colon \bp(V^*)\longrightarrow 
\bp ( (\sym^n V^*)^G )$$
defined by $\varphi( [f] ) = [\prod_{\sigma\in G} \sigma \cdot f ] $ 
is nondegenerate
for every finite-dimensional complex representation $V$ of $G$
 if and only if $n $ is either a prime number or $n=4$.}

\vspace{1pc}
We first prove the sufficiency of the condition.
Let $n=p$, a prime number and
$\zeta$ be a primitive  $p$th root of unity. Let $\{X_{ij}\mid 
i=1,\ldots, m_j, j=0, \ldots, p-1\}$ be a basis of $V^*$ such that
 $\sigma X_{ij}=\zeta^jX_{ij}\hbox{ for all } i,j$. 
 We prove that for any $p$-tuple $\a =(a_0, \ldots, a_{p-1})$ of
nonnegative integers 
such that $\sum_{i=0}^{p-1}a_i=p$  and $\sum_{i=0}^{p-1} a_i i 
\equiv 0 \pmod p$, and for any subset $A_j$ of $\{1,\ldots, m_j\}$ 
of cardinality $a_j$, the monomial $\prod_{j=0}^{p-1} \prod_{i\in A_j} X_{ij} \prod_{j=0}^{p-1}\prod_{i\in A_j} Y_{ij}$ occurs with 
nonzero coefficient in 
$$   \prod_{t=0}^{p-1} \sigma^t\left( \sum_{j=0}^{p-1}\sum_{i=1}^{m_j} X_{ij}Y_{ij}  \right)  $$

The argument depends only on weight computation. So, we may assume 
that $X_{ij}=X_{kj}\ \hbox{ for all } i,j,k$.

Let $\pi:{\bz}\to {\bz}/p\bz$ be the natural map. 
Let $x\in {\bz}/p{\bz}$. 
Define $b_r =\sum_{i=0}^r a_i$.
Define 
$S_{x,\a}=\{\sigma\in S_p\mid \pi \left(\sum_{r=0}^{p-1}
\left(\sum_{i=1+b_{r-1}}^{b_r} \sigma(i)\right)r\right)=x\}$.  
Then, we have 
\begin{lem}
 For any $x,y\in {\bz}/p{\bz}$ 
such that $y=u\cdot x$ for some nonzero $u$ in ${\mathbf Z}/p{\mathbf Z}$, 
 $(1)\  \left| S_{x,\a}\right| =\left| S_{y,\a}\right|$, and 
$ (2)\  \left| S_{\pi(0),\a}\right|\neq \left| S_{\pi(1),\a }\right|$.
\end{lem}

\noindent{\bf Proof of }(1) The map 
$\tilde{u}:S_{x,\a}\to S_{y,\a}$ given by $ \tilde{u}(\sigma)(i)=u\pi(\sigma(i)), i=0, \ldots p-1$ gives a bijection.

\noindent{\bf Proof of} (2). We first observe that $$p  
\hbox{ divides } 
\left| S_{x,\a} \right|\ \hbox{for all }x\in {\mathbf Z}/p{\mathbf Z} 
\hbox{ and } \forall \a\eqno(*)$$
For a proof, let $z\in {\mathbf Z}/p{\mathbf Z}$ and let 
$\sigma\in S_{x,\a}$. Then the permutation given by  
$(\tilde{z}\circ \sigma)(\pi(i))=\sigma(\pi(i))+z$
belongs to $S_{x,\a}$ as $\sum_{r=0}^{n-1}a_rr \equiv 0 \pmod p. $ 
Thus $p \hbox{ divides }\left| S_{x,\a}  \right|\ \forall x \mbox{ and } \a$. 
On the other hand, $$ \bigcup_{x\in {\mathbf Z}/p{\mathbf Z}}S_{x,\a}
=S_p \eqno(**)$$
From $(*)$ and $(**)$, we see that if $ \left| S_{\pi(0),\a}\right|
=\left| S_{\pi(1),\a} \right|$, then $p\cdot \left| S_{0,\a}\right|
=p\mbox{!}$.
 Hence $p$ divides $(p-1)$! by $(*)$ which is a contradiction for a prime $p$.

\begin{cor}
Let $p$ be prime. For any $p$-tuple of
nonnegative integers   
$\a=(a_0,\ldots, a_{p-1})$  such that $\sum_{r=0}^{p-1}a_r
=p$ and $\sum_{r=0}^{p-1}a_rr\equiv 0 \pmod p$; the polynomial 
\break 
$\sum_{\sigma\in S_p} T^{\sum_{r=0}^{p-1}\left( \sum_{i=1+ b_{r-1}}^{b_r} \sigma(i)\right)r}$ is not divisible by the $p$th 
cyclotomic polynomial.
\end{cor}

Thus, we have shown that the coefficient of $ X^\a Y^\a$ in 
$\prod_{i=0}^{p-1}\sigma^i \left( \sum_{r=0}^{p-1}X_{r}Y_r\right)$  
is non-zero for all $\a$.

The nondegeneracy for both groups of order 4 can be checked
explicitly.

\vspace{1pc}
The necessity of the condition in the Theorem  is a consequence of 
the following: 
\begin{prop}
 Let $G$ be a finite  abelian group 
whose order $n$ is a composite number  different from $4$. Let $V$ be the regular representation of $G$ over \bc. Then the morphism
$$ \varphi\colon \bp(V^*)\longrightarrow \bp ( (\sym^n V^*)^G )$$
defined by $\varphi( [f] ) = [\prod_{\sigma\in G} \sigma \cdot f ] $ is
degenerate.
\end{prop}
\vspace{6pt}
{\bf Proof:}\quad
This proposition will be proved by considering three cases:
 cyclic groups, groups which are direct products of at least
3 copies of the group of order 2, and the other finite abelian groups.
In all the cases we exhibit a monomial $M$ which 
is not in the span of  image of $\varphi$, proving degeneracy 
of the map. 

\noindent{\bf Case I}: $G$ is a cyclic group of order $n$.
There are two subcases to be considered here.
Let  $p$ denote the least prime divisor of $n$, and let
$q = n/p$.

For construction of $M$, we use a basis $\{X_\lambda \mid
\lambda \in G \}$ of common eigen vectors of $G$. Note that we abuse 
notation, by using elements of $G$, instead of the dual  $\hat G$
for indexing this basis.

\textbf{Subcase Ia.}\quad  $p|q$. Take
$M= X_0X_1^{q-2}X_{n-q+2}X_q^{(p-1)q}$

\textbf{Subcase Ib}.\quad $\gcd(p,q) =1$. Take
$M= X_0^{q-2}X_{q-p}X_{n-q+p} X_p^{(p-1)q}$

\vspace{6pt}
For the other cases we decompose the group using structure theorem
of finite abelain groups.
$G = \bigoplus_{i=1}^r \bz/q_i\bz$, with 
$q_i|q_{i+1},\ i=1,2,\ldots r-1$.

\vspace{6pt}
\noindent{\bf Case II}. $r\geq2$ and $q_r\geq3$

\vspace{3pt}
Let $H = \bigoplus_{i=2}^r \bz/q_i\bz$ have index $q$ in $G$.
That is $n  = q q_1$. 
Take $M = X^{q-2}_{0,0,\ldots,0}X_{0,1,1,\ldots,1}X_{0,q_2-1,\ldots,q_r-1}X^{q(q_1-1)}_{1,1,\ldots,1}$.

\vspace{6pt}
\noindent{\bf Case III}
 $G$ is a product of 3 or more copies of $\bz/2\bz$.
Now we express $G$ as $\bz/2\bz \oplus U$, with $U$ a subgroup
of order $q = n/2$.
Take $M= (\prod_{v\in U} X_{0,v})X_{1,1,\ldots,1}^q$.

Now we give the arguments for each case.

\textbf{Case Ia}.
Then the monomial  is  $G$-invariant  of degree $n$, and  it  can be expressed in only
 the following two ways as product of $p$ monomials each of which is a
 $H$-invariant monomial of degree $q$,  where $H$ is the unique
subgroup of order $q$ in $G$.

\begin{itemize}
\item[(i)] $M =\prod_{i=1}^p M_i$ where 
 $M_1=X_0X_1^{q-2}X_{n-q+2};\ M_i=X_q^q,\ i=2,\ldots, p$.

\item[(ii)] $M=\prod_{i=1}^p M_i'$ where  $M_1'=X_qX_1^{q-2}X_{n-q+2};\ M_2'=X_0X_q^{q-1}$

$M_i'=X_q^q, \ i=3,4,\ldots, p$.
\end{itemize}
It is clear that $X_1^{q-2}$ and $X_{n-q+2}$ must be a factor of 
a common $H$-invariant monomial of degree $q$. So either such 
a monomial must be  $X_0X_1^{q-2}X_{n-q+2}$ or $X_qX_1^{q-2}X_{n-q+2}$.

Now we prove that the coefficient of $M_i(X)M_i(Y)$ in $\prod_{j=0}^{q-1}\sigma^{pj}\left(\sum_{r=0}^{n-1} X_r Y_r\right)$ is the same as the 
coefficient of ${M_i'(X) M_i'(Y)}$ in $\prod_{j=0}^{q-1} \sigma^{pj} \left(\sum_{r=0}^{n-1} X_r Y_r\right)$.

We prove that the monomial  
 $X_0X_1^{q-2}X_{n-q+2}X_q^{(p-1)q} Y_0Y_1^{q-2}Y_{n-q+2}Y_2^{(p-1)q}$
does not appear in the product  
$\prod_{i=0}^{p-1}\sigma^i \left(\prod_{j=0}^{q-1}\sigma^{pj}\left(\sum_{r=0}^{n-1}X_rY_r\right)\right)$.
 
Let $\zeta$ be a primitive $n$th 
root of 1. Then, the coefficient of ${M_1(X) M_1(Y)}$ in 
$\displaystyle {\prod_{j=0}^{q-1} \sigma ^{pj}\left( \sum_{r=0}^{n-1} X_r Y_r\right)}$  is 
\begin{eqnarray*}
&&=\sum_{\sigma\in S_q}\zeta^{p\left(\sigma(0)\cdot 0+\left(\sum_{t=1}^{q-2}\sigma(t)\right)\cdot 1+\sigma(q-1)(n-q+2)\right)  }\\
&&=\sum_{\sigma\in S_q}\zeta^{p\left(\sigma(0)\cdot q+\left(\sum_{t=1}^{q-2}\sigma(t)\right)\cdot 1+\sigma(q-1)(n-q+2)\right)  }\\
&&\left( \mbox{because } \zeta^{jpq}=1 \forall j\in \{0,1,2,\ldots,q-1\}\right)\\
&&=\mbox{Coefficient of } M_1'(X)M_1'(Y)\mbox{ in } \prod_{j=0}^{q-1}\sigma^{pj}\left(\sum_{r=0}^{n-1}X_r Y_r\right)
\end{eqnarray*}
Similarly, the coefficient of $M_2(X)M_2(Y)= $ the coefficient of $M_2'(X)M_2'(Y)$.

Since $M_i=M_i'$ for $i=3,4,\ldots, p$, the coefficients of 
 $M_i(X)M_i(Y)$  and $M_i'(X)M_i'(Y)$ in  
$\prod_{j=0}^{q-1}\sigma^{pj} \left(\sum_{r=0}^{n-1}X_rY_r\right)$
are equal.

Let $c_i$ be the coefficient of $M_i(X)M_i(Y)$ in $\prod_{j=0}^{q-1}\sigma^{ pj}\left(\sum_{r=0}^{n-1} X_r Y_r\right)$. Let $c_i'$ be the 
coefficient of $M_i'(X)M_i'(Y)$. Then, $c_i=c_i'$ by 1. Then, the 
coefficient of  $X_0X_1^{q-2}X_{n-q+2}X_q^{(p-1)q}$ in 
$\prod_{i=0}^{p-1}\sigma^i\left(\prod_{j=0}^{q-1}\sigma^{pj}
\left(\sum_{r=0}^{n-1}X_rY_r\right)\right)$ 
is 

\noindent $\Big(\sum_{i=0}^{p-1}\zeta^{in+q^2\sum_{k\neq i}k}
+\sum_{i,j=0\atop i\neq j}^{p-1} \zeta^{i(n+q)+j(q(q-1))
 + q^2 \sum_{k=0\atop k\neq i,j }^{p-1} k}\Big)c$ 
where $c=\prod_{i=1}^p c_i$.

\vspace{2pt}
But \begin{eqnarray*}&&\sum_{i=0}^{p-1}
 \zeta^{in+q^2\sum_{k\neq i}k}
    +\sum_{i,j=0\atop i\neq j}^{p-1} \zeta^{i(n+q)+jq(q-1)
  + q^2 \sum_{k=0\atop k\neq i,j }^{p-1} k}\\
&&{}=p+\sum_{\textstyle {i\neq j\atop i,j=0} }^{p-1}\zeta^{(i-j)q}
\quad(\hbox{since } \zeta^n=\zeta^{q^2}=1)\\
&&{}=p+p\left(\sum_{t=1}^{p-1}\zeta ^{tq}\right)\\
&&{}=p\left( \sum_{i=0}^{p-1}\zeta^{ tq}  \right)=0.
\end{eqnarray*}
\noindent {\bf Case Ib}.

By a similar argument, one can show that if $\gcd(p,q)=1$,
 the coefficient of 
$X^{q-2}_0X_{q-p}X_{n-q+p}X_p^{(p-1)q }
Y^{q-2}_0Y_{q-p}Y_{n-q+p}Y_p^{(p-1)q }$ in 
$\prod_{i=0}^{n-1}\sigma^i \left(\sum_{r=0}^{n-1}X_rY_r\right)$ 
is zero.

\vspace{1pc}
\noindent{\bf Case II.}

Let $(M_1,\ldots, M_{q_1})$ be a $q_1$-tuple of monomials each
of which is $G$-invariant such that their product $\prod _{i=1}^{q_1} M_i$ is $M$. There is a unique $i_0$ and a unique $j_0$ such that
the variables $X_{0,q_2-1,\ldots,q_r-1}$  and 
$X_{0,1,1,\ldots,1}$  occur respectively  in the monomials 
$M_{i_0}$ and and $M_{j_0}$. As we are in the case
$q_r \geq3$, these two variables are distinct.
So for each $t$ the monomial $M_t$ would be $G$-invariant
if and only if $i_0 =j_0$.

For a fixed $i$ consider the set of $q_1$-tuples of monomials defined by
$$S_i = \{ (M_1,\ldots M_{q_1}) \mid X_{0,q_2-1,\ldots,q_r-1} 
\hbox{ occurs in } M_i\}$$
and fixing a $j$ define the subset
$$S_{i,j} = \{ (M_1,\ldots M_{q_1})\in S_i \mid X_{0,1,1,\ldots,1} 
\hbox{ occurs in } M_j\}$$

Let $(M_1,\ldots M_{q_1})\in S_{i,i}$. Let $j\neq i$.
Define $$ \begin{array}{rlll}M_t' &=&  M_t &  
 \hbox{for all } t\neq i,j\\[2pt]
  M_i'&=&\displaystyle
  \frac{M_i}{X_{0,1,1\ldots,1}} X_{1,1,\ldots,1}\\[11pt]
  M_j'&=& \displaystyle
  \frac{M_j}{X_{1,1,1\ldots,1}} X_{0,1,\ldots,1}.\\
\end{array}$$
Then the map $S_{i,i}\to S_{i,j}$
sending $(M_1,\ldots,M_{q_1})$ to $(M_1',\ldots,M_{q_1}')$ is 
bijective.
If $a_i$ is  the  sum of coefficients of $M(X)M(Y)$ in 
$\prod_{i=0}^{q_1-1} \sigma^i \prod _{\tau\in H}
\tau (\sum _{\lambda\in \hat G} X_\lambda Y_\lambda)$
contributed by the  elements 
$ (N_1,N_2,\ldots,N_{q_1} )\in S_{i,i} $ such that $N_i = M_i$,
then the sum of coefficients of $M(X)M(Y)$ contributed by the larger
set
$$\{(N_1,N_2,\ldots,N_{q_1} )\in S_i \mid N_i = M_i \}$$ is 
$a_i (\sum _{j=0}^{q_1-1} \zeta ^{i-j})$, with $\zeta$ being a primitive $q_1$-th root of unity. But this sum is zero as it is the sum of \textit{all} the $q_1$-th roots of unity.

This proves that 
in the product 
$\prod_{i=0}^{q_1-1} \sigma^i \prod _{\tau\in H} 
\tau (\sum _{\lambda\in \hat G} X_\lambda Y_\lambda)$
the coefficient of $M(X)M(Y) $   
is zero.

\vspace{6pt}
\noindent{\bf Case III}
 
The proof proceeds along the same lines as in Case II, making use
of the Lemma 2.

\end{document}